\documentclass[a4paper,10pt]{amsart}
\usepackage{amsmath,amsfonts,amssymb, amsthm, xypic}
\usepackage{epsfig,psfrag}
\usepackage[all]{xy}

\def\d{\mathbf{d}}
\def\i{\mathbf{i}}

\def\Q{\mathbb{Q}}

\def\N{\mathbb{N}}

\def\P{\mathbb{P}}

\def\mQ{\mathcal{Q}}
\def\L{\mathcal{L}}

\newtheorem*{lem}{Lemma}

\newtheorem*{cor}{Corollary}
\newtheorem*{prop}{Proposition}

\begin{document}

\title[Addendum]{Addendum to ``Canonical bases for quantum generalized Kac-Moody algebras''}

\author{Seok-Jin Kang$^*$ and Olivier Schiffmann$\dag$}
\date{}

\address{$^*$ Department of mathematical sciences and \\
research institute for mathematics, Seoul national university,\\
San 56-1 Shinrim-dong, Kwanak-ku\\
Seoul 151-747, Korea}
\email{sjkang@kias.re.kr}
\address{$^\dag$ DMA, \'Ecole Normale Sup\'erieure, 45 rue d'Ulm,
          75230 Paris Cedex 05-FRANCE}
\email{schiffma@dma.ens.fr}

\begin{abstract}

We provide some necessary details to several arguments appearing in our previous paper
 ``Canonical bases for quantum generalized Kac-Moody algebras''. 
\end{abstract}

\maketitle

\paragraph{} As pointed out to us by several people, some of the arguments in our paper \cite{KS} are too sketchy and at places incomplete. The first purpose of this addendum is to fill in the missing details and provide complete proofs. The second purpose is to explain the precise relation between \cite{KS} and other preexisting works on canonical bases for quivers with loops (mainly the paper \cite{LusTight}). 

\vspace{.1in}

We will freely make use of the notation from \cite{KS}. In particular, we have fixed a quiver $Q$, denoted by $I$ its set of vertices, by $I^{im}$ the subset of $I$ consisting of vertices with edge loops (\textit{imaginary vertices}) and by $I^{re}$ the complement of $I^{im}$ (\textit{real vertices}). Let $\Omega$ stand for an orientation of $Q$. We denote by $r_{ij}$ the number of oriented arrows going from $i$ to $j$ and by $c_i$ the number of edge loops at the vertex $i$.
For any dimension vector $\mathbf{d} \in \N^I$ there is a representation space $E_{\mathbf{d}}$ for the quiver $Q$ over the algebraic closure $\overline{\mathbb{F}_q}$ of the finite field $\mathbb{F}_q$. The reductive group $G_{\mathbf{d}}=\prod_i GL(d_i, \overline{\mathbb{F}_q})$ acts on $E_{\mathbf{d}}$, and the orbits are in bijection with the set of all (nilpotent) representations of $Q$ of dimension $\mathbf{d}$. To any sequence $\mathbf{i}=(i_1, \ldots, i_r)$ of vertices of $Q$ are attached two kinds of $I$-graded flags varieties~:

$$\mathcal{F}_{\mathbf{i}}=\{D_{\bullet}\;|0 \subseteq D_1 \subseteq
\cdots \subseteq D_r=V_\d;\; \mathbf{dim}
(D_l/D_{l-1})=\epsilon_{i_l}\},$$

and

$$\mathcal{F}^{im}_{\mathbf{i}}=\{D_{\bullet}\;|0 \subseteq D_1 \subseteq
\cdots \subseteq D_s=V_{\mathbf{d}^{im}};\; \mathbf{dim}
(D_l/D_{l-1})=\epsilon_{i^{im}_l}\},$$

Here $\mathbf{d} \in \N^I$ is defined by $d_j=\#\{l\;|\; i_l=j\}$ and $V_{\mathbf{d}}$ is some fixed $I$-graded $\overline{\mathbb{F}_q}$-vector space of dimension $\mathbf{d}$. Similarly, $\mathbf{i}^{im}$ is the sequence of imaginary vertices obtained from $\mathbf{i}$ by removing all real vertices, $\mathbf{d}^{im}$ is the projection of $\mathbf{d}$ to $\N^{I^{im}}$, and $V_{\mathbf{d}^{im}}$ is the direct summand of $V_{\mathbf{d}}$ corresponding to imaginary vertices.

\vspace{.1in}

We also consider the corresponding incidence varieties

$$\widetilde{\mathcal{F}}_{\mathbf{i}}=\{(x,D_{\bullet})\;|x(D_i)
\subseteq D_{i-1}\} \subset E_{\mathbf{d}} \times
\mathcal{F}_{\mathbf{i}},$$
$$\widetilde{\mathcal{F}}^{im}_{\mathbf{i}}=\{(x,D_{\bullet})\;|x(D_i)
\subseteq D_{i-1} \oplus \bigoplus_{i \in I^{re}} V_i\}
\subset E_{\mathbf{d}} \times \mathcal{F}_{\mathbf{i}}^{im}.$$

\vspace{.1in}

There is a sequence of proper and $G_{\mathbf{d}}$-equivariant maps

$$
\xymatrix{
\widetilde{\mathcal{F}}_{\mathbf{i}}
\ar[r]^{\pi_1'} &
\widetilde{\mathcal{F}}_{\mathbf{i}}^{im}  
\ar[r]^{\pi_1}&
E_{\mathbf{d}} }
$$

By the Decomposition Theorem \cite{BBD}, the complex $\pi'_{1!}((\overline{\mathbb{Q}}_l)_{\widetilde{\mathcal{F}}_{\mathbf{i}}})[\mathrm{dim}\;\widetilde{\mathcal{F}}_{\mathbf{i}}]$ is semisimple. We let
$\mathcal{T}_{\mathbf{i}}$ be the set of all simple perverse sheaves
appearing (possibly with a shift) as simple consituents of $\pi'_{1!}((\overline{\mathbb{Q}}_l)_{\widetilde{\mathcal{F}}_{\mathbf{i}}})[\mathrm{dim}\;\widetilde{\mathcal{F}}_{\mathbf{i}}]$.

\vspace{.2in}

\paragraph{\textbf{1. Perversity statement in Proposition~4.1.}} Proposition~4.1. of \cite{KS} states that for any $\mathbb{P} \in \mathcal{T}_{\mathbf{i}}$ the semisimple complex $\pi_{1!}(\mathbb{P})$ is in fact perverse. The proof, as written there, is based on a ``standard argument''. This is at best misleading. 
Before going any further, let us recall a few of the notations. Let $E'_{\d}\subset\bigoplus_{i \neq j} \mathrm{Hom}(V_i,V_j)^{r_{ij}}$
and $E_{\d}''\subset\bigoplus_{i} \mathrm{Hom}(V_i,V_i)^{c_i}$ be the
set of nilpotent representations and let
$u: E_{\d} \to E'_{\d}$, $t: E_{\d} \to E''_{\d}$ be the projections.
Set $\mathcal{G}_{\i}= (u \times
Id)(\widetilde{\mathcal{F}}_{\i}) \subset E'_{\d} \times
\mathcal{F}_{\i}$ and
$\mathcal{G}^{im}_{\i}= (u \times Id)(\widetilde{\mathcal{F}}^{im}_{\i})
\subset E'_{\d} \times \mathcal{F}^{im}_{\i}$, so that there is a
commutative  diagram
$$\xymatrix{
\widetilde{\mathcal{F}}_{\mathbf{i}}
\ar[r]^{\pi_1'} \ar[d]_{u \times Id}&
\widetilde{\mathcal{F}}_{\mathbf{i}}^{im}
\ar[d]^{u\times Id} \ar[r]^{\pi_1}&
E_{\mathbf{d}} \\
\mathcal{G}_{\mathbf{i}} \ar[r]^{s} &
\mathcal{G}_{\mathbf{i}}^{im}
}$$
Observe that the vertical maps are vector bundles, and that
$u \times Id:\widetilde{\mathcal{F}}_{\mathbf{i}} \to \mathcal{G}_{\i}$
is the pullback by $s$ of the bundle $u\times Id:
\widetilde{\mathcal{F}}_{\mathbf{i}}^{im}
\to \mathcal{G}_{\i}^{im}$. Hence
$\pi'_{1!}((\overline{\mathbb{Q}}_l)_{\widetilde{\mathcal{F}}_{\i}})=
\pi_{1!}'(u\times Id)^*((\overline{\mathbb{Q}}_l)_{\mathcal{G}_{\i}})=
(u\times Id)^*s_!((\overline{\mathbb{Q}}_l)_{\mathcal{G}_{\i}})$.  In
particular, any
of the simple perverse sheaves in $\mathcal{T}_\i$ is of the form
$IC(X,\mathfrak{L})$ with $X = (u\times Id)^{-1}(Y)$ for a smooth
irreducible  subvariety $Y \subset \mathcal{G}_{\i}^{im}$ and
$\mathfrak{L}=
(u \times Id)^* \mathfrak{K}$ for an irreducible local system
$\mathfrak{K}$ on $Y$. From that it is deduced in \cite{KS} that the restriction of 
$\pi_1$ to $\overline{X}$ is semismall, and therefore that, ``by standard arguments'',
$\pi_{1!}(\mathbb{P})$ is perverse. When $\mathbb{P}$ is a local system, this is certainly a standard fact; but this is not necessarily the case here. Nevertheless, perversity of $\pi_{1!}(\mathbb{P})$ follows from the fact that the restriction of $\pi_{1}$ to \textit{each strata} of $\mathbb{P}$ is semismall. To see this,
let $\overline{Y}=Y_1 \sqcup Y_2 \cdots \sqcup Y_n$ be a decomposition of $\overline{Y}$ into (smooth irreducible) strata for $IC(Y,\mathcal{K})$. Then $X_i:=(u \times Id)^{-1}(Y_i)$ for $i=1, \ldots, n$ are the strata for $IC(X,\L)$ (i.e. the restriction of $H^j(IC(X,\L))$ to each $X_i$ is a local system).

\vspace{.1in}

By the same argument as for the whole of $\overline{X}$, the restricted projection map
$(\pi_1)_{|X_i}: X_i \to \pi_1(X_i)$ is semismall for every $i$ (and in particular, $dim\;(X_i)=dim\;(\pi_1(X_i))$).

\vspace{.1in}

Set $e_i=codim_{\overline{X}} X_i=codim_{\pi_1(\overline{X})} \pi_1(X_i)$ and $e_0=dim\;X$. We have
$$H^{>e_i}(IC(X,\L)[e_0]) \equiv 0\qquad \text{on}\;X_i$$
for all $i$. By semismallness of $(\pi_1)_{|X_i}$ and since $H^j(IC(X,\L)[e_0])$ restricts to a local system on each $X_i$, we have for every $j=0, \ldots, e_i$
$$codim_{\pi_1(X_i)} \;supp\; \big( H^j(\pi_{1!}(IC(X,\L)[e_0]_{|X_i}\big) \geq j-e_i$$
from which it follows that
$$codim_{\pi_1(\overline{X})} \;supp\; \big( H^j(\pi_{1!}(IC(X,\L)[e_0]_{|X_i}\big) \geq j.$$

\vspace{.1in}

By the long exact sequence in $H^*$ associated to $\pi_{1!}$ and the triangles induced by successive closed inclusions of strata we have
$$supp \big( H^j(\pi_{1!}(IC(X,\L)[e_0]))\big)=\bigcup_isupp \big( H^j(\pi_{1!}(IC(X,\L)[e_0]_{|X_i}))\big).$$
Hence $\pi_{1!}(IC(X,\L))$ satisfies the support condition for a perverse sheaf. The Verdier dual
$\mathcal{D}IC(X,\L)$ is a perverse sheaf owhich also belongs to $\mathcal{T}_{\mathbf{i}}$ and the same argument may be applied to it. Hence it also satisfies the support condition, and therefore $\pi_{1!}(IC(X,\L))$ is indeed perverse.

\vspace{.15in}

\noindent
\textbf{Remark.} A semismallness result very close to the one given in \cite{KS}, Proposition~4.1. appears in the paper \cite{LusTight}.

\vspace{.2in}

Let us state the following corollary of Proposition~4.4., which is used (without due statement) in the course of the proof of Proposition~5.1.~:

\begin{cor} For any $\mathbf{i}$ and any $R \in \mathcal{T}_{\mathbf{i}}$ we have
$$supp(R) \subset \pi^{-1}_{1} (supp ( \pi_{1!}(R))).$$
\end{cor}

\noindent
\textit{Proof.} This comes from the fact that the map $\pi_1$ is semismall. More precisely, each $R$ as above is of the form $IC(X,\mathcal{L})$ where $X=(u^{-1} \times Id) (Y)$ and $\mathcal{L}=(u \times Id)^* \mathcal{K}$ for some subvariety $Y$ of $\mathcal{G}_{\mathbf{i}}^{im}$ and local system $\mathcal{K}$ on $Y$. Let $X^{reg}$ be the open subset of $X$ consisting of points $(x,D_{\bullet})$ for which all the maps in $x$ associated to edge loops are regular nilpotent. Then
$$\overline{X^{reg}}=\overline{X}, \qquad \overline{\pi(X^{reg})}=\overline{\pi_1(X)}=\overline{\pi_1(\overline{X})}.$$
But $(\pi_1)_{|X^{reg}} : X^{reg} \to \pi_1(X^{reg})$ is an isomorphism, hence
$$supp(\pi_{1!}(R)) \supset \overline{ \pi_1(X^{reg})}=\overline{\pi_1(\overline{X})}=\overline{\pi_1(Supp(R))}.$$

\vspace{.2in}

\paragraph{\textbf{2. Fourier-Deligne transform.}}
Recall that for each $\mathbf{i}$ we have defined the set of semisimple perverse sheaves $\mathcal{P}_{\mathbf{i}}=\{\pi_{1!}(P)\;|P \in \mathcal{T}_\mathbf{i}\}$. We let $\mathcal{P}_{\d}=
\bigcup_{\mathbf{i}} \mathcal{P}_{\mathbf{i}}$ where the sum ranges over
all  sequences $\mathbf{i}$ such that $\sum_l \epsilon_{i_l} =\d$.
Finally we denote by $\mathcal{Q}_{\d}$ the category of complexes which
are direct sums of shifts of elements in $\mathcal{P}_{\d}$. Thus
$\mathcal{Q}_{\d}$ is a full subcategory of
$\mathcal{Q}_{G_{\d}}(E_{\d})$. 

\vspace{.1in}

Let $\Omega_1, \Omega_2$ be two orientations of our quiver $Q$. We may write $\Omega_1=\Omega_0 \sqcup \Omega'$ and $\Omega_2=\Omega_0 \sqcup \overline{\Omega'}$, and assume that all edge loops belong to $\Omega_0$. To stress the choice of an orientation, we add in all notations a subscript $\Omega$, as in $\mQ_{\d,\Omega}, \mathcal{P}_{\d,\Omega}$, etc.
Let $\overline{\mQ}_{\d,\Omega_i}$, for $i=1,2$ stand for the category generated by all \textit{simple} constituents of the $\pi_{\mathbf{i} !}(\mathbf{1}_{\widetilde{\mathcal{F}}_{\mathbf{i}}})$. Thus the semisimple perverse sheaves in $\mQ_{\d,\Omega_i}$ are certain direct sums of objects in $\overline{\mQ}_{\d,\Omega_i}$. The Fourier-Deligne transform of Lusztig (see \cite{Lusbook}, Chap. 11) defines an equivalence of categories ${\Phi}_E: \overline{\mQ}_{\d,\Omega_1} \stackrel{\sim}{\to} \overline{\mQ}_{\d,\Omega_2}$. Note that we only perform the Fourier-Deligne transform on arows which are not edge loops so the aguments in \cite{Lusbook} apply \textit{verbatim}. In the course of the proof of Theorem~4.1. in \cite{KS} we have used without due justification the following~:

\begin{prop} The equivalence $\Phi_E$ restricts to an equivalence
$$\Phi_E: \mQ_{\d,\Omega_1} \stackrel{\sim}{\to} \mQ_{\d,\Omega_2}.$$
\end{prop}

\vspace{.1in}

\noindent
\textit{Proof.} Recall that any object in $\mathcal{P}_{\d,\Omega_i}$ is obtained as $\P=\pi_{1!}(P)$, where
$P$ is a simple perverse sheaf on $\widetilde{\mathcal{F}}_{\i}^{im}$ and $\pi_1: \widetilde{\mathcal{F}}_{\i}^{im} \to E_{\d}$ is the projection. Moreover, any occuring such simple perverse sheaf $P$ is the pullback under the vector bundle $\widetilde{\mathcal{F}}_{\i}^{im} \to \mathcal{G}_{\i}^{im}$ of some simple perverse sheaf on $\mathcal{G}_{\i}^{im}$. The basic idea is to define several Fourier transforms, for $E_{\d}$, $\widetilde{\mathcal{F}}_{\i}^{im}$, $\mathcal{G}_{\i}^{im}$, and to show that they are all compatible in some sense. However, it is better for tchnical purposes to use spaces which are slightly larger than $\widetilde{\mathcal{F}}_{\i}^{im}$ and $\mathcal{G}_{\i}^{im}$.

\vspace{.1in}

Let us write an element $x \in E_{\d}$ as $x=(x',x^{nil})$, where $x$ is the direct sum of all maps which do not correspond to edge loops, and $x^{nil}$ is the direct sum of all maps corresponding to edge loops. Accordingly, there is a decomposition $E_{\d}=E'_{\d} \times E^{nil}_{\d}$. For $\i$ a sequence of vertices, put
$$\mathcal{H}^{im}_{\i}=\{(x',x^{nil},D_{\bullet})\;|\; x^{nil}(D_l) \subset D_{l-1}\;for\;all\;l\} \subset E_{\d} \times \mathcal{F}^{im}_{\i}.$$
The difference between $\widetilde{\mathcal{F}}_{\i}^{im}$ and $\mathcal{H}^{im}_{\i}$ is that in $\widetilde{\mathcal{F}}_{\i}^{im}$ we impose in addition that $x_h(D_l) \subset D_{l-1}$ when $x_h$ is an arrow joining two \textit{distinct} imaginary vertices. In particular, $\widetilde{\mathcal{F}}_{\i}^{im}$ is a closed subset of $\mathcal{H}^{im}_{\i}$.
Finally, we put 
$$ \mathcal{K}^{im}_{\i}=E'_{\d} \times \mathcal{F}^{im}_{\i},$$
$$E_{\d}^{reg}=E'_{\d} \times E^{reg, nil}_{\d}$$
where $E^{reg,nil}_{\d} \subset E^{nil}_{\d}$ stands for the open subset of regular nilpotent elements (i.e. the $x_h$ are regular for all the edge loops $h \in \Omega$). Again, the difference between $\mathcal{K}^{im}_{\i}$ and $\mathcal{G}^{im}_{\i}$ is that in $\mathcal{G}^{im}_{\i}$ one in addition requires the maps corresponding to arrows linking two \textit{distinct} imaginary vertices to be compatibe with the flag.

\vspace{.1in}

All the above varieties are related by the following diagram~:

\begin{equation}\label{Diagram}
\xymatrix{
\mathcal{H}^{im}_{\i} \ar[r]^{\pi_1} \ar[d]^{e} & E_{\d}\\
\mathcal{K}^{im}_{\i} & E_{\d}^{reg} \ar[u]_j \ar[l]^{\rho}
}
\end{equation}

\noindent
where $\pi_1$ is the (proper) projection, $j$ is the open embedding, $\rho$ is the smooth map which associates to a regular element $x=(x',x^{nil})$ the pair $(x',D_{\bullet})$ where $D_{\bullet}$ is the only flag of type $\i^{im}$ compatible with $x^{nil}$, and finally $e$ is the natural vector bundle map.

\vspace{.1in}

Now let as above $\Omega_1, \Omega_2$ be two different orientations. Then we have vector bundle maps
$$\mathcal{H}^{im}_{\i,\Omega_i} \to \mathcal{H}^{im}_{\i,\Omega_0},$$
$$ \mathcal{K}^{im}_{\i,\Omega_i} \to \mathcal{K}^{im}_{\i,\Omega_0},$$
$$E_{\d,\Omega_i} \to E_{\d,\Omega_0},$$
$$ E^{reg}_{\d,\Omega_i} \to E^{reg}_{\d,\Omega_0}$$
for $i=1,2$. All squares of the form
$$\xymatrix{
X_{\Omega_1} \ar[r] \ar[d] & X_{\Omega_0} \ar[d]\\
Y_{\Omega_1} \ar[r] & Y_{\Omega_0}
}
$$ 
built from the diagram (\ref{Diagram}) commute and are cartesian. Moreover, there are natural pairings
$$X_{\Omega_1} \times_{X_{\Omega_0}} X_{\Omega_2} \to k$$
and these are all compatible with the maps $X_{\Omega_i} \to Y_{\Omega_i}$ in (\ref{Diagram}).

The Fourier transform defines an equivalence of categories
$$\Phi_X: D^b(X_{\Omega_1}) \stackrel{\sim}{\to} D^b(X_{\Omega_2})$$
for $X \in \{\mathcal{H},\mathcal{K},E,E^{reg}\}$.

\vspace{.1in}

\begin{lem} We have
\begin{enumerate}
\item[i)] $j^* \circ \Phi_E \simeq \Phi_{E^{reg}} \circ j^*$,
\item[ii)] $\rho^* \circ \Phi_{E^{reg}} \simeq \Phi_{\mathcal{K}} \circ \rho^*$,
\item[iii)] $e^* \circ \Phi_{\mathcal{K}} \simeq \Phi_{\mathcal{H}} \circ e^*$,
\item[iv)] $\pi_{1!} \circ \Phi_{\mathcal{H}} \simeq \Phi_{{E}} \circ \pi_{1!}.$
\end{enumerate}
\end{lem}
\noindent
\textit{Proof.} For i), it is the compatibility of the Fourier transform with the restriction to an open subset of the base; for ii) and iii) it is the compatibility of the Fourier transform under smooth pullback of the base.
These are all special case of the compatibility of Fourier-Deligne transform with base change, see e.g. \cite{Laumon}.

We prove iv).  For simplicity, we drop all unnecessary indices (such as $\d$, $\i$, etc..). Let $\P$ be an object of $D^b(\mathcal{H}_{\Omega_1})$, and $\Q$ an object of $D^b(E_{\Omega_1})$. By definition, 
$$\Phi_{\mathcal{H}}(\P)=p_{2!}( p_{1}^*(\P) \otimes \mathcal{L}_{\mathcal{H}}), \qquad
\Phi_{E}(\Q)=p'_{2!}( {p'}_{1}^*(\Q) \otimes \mathcal{L}_{E}),$$ 
where $p_1, p_2$ and $p'_1, p'_2$ are the relevant vector bundle maps. We have

\begin{equation*}
\begin{split}
\pi_{1!} \Phi_{\mathcal{H}}(\P)&=\pi_{1!}p_{2!}(p_1^*(\P) \otimes \mathcal{L}_{\mathcal{H}})\\
&={p'}_{2!}\pi_{1!}(p_1^*(\P) \otimes  \mathcal{L}_{\mathcal{H}}).
\end{split}
\end{equation*}
Note that $\mathcal{L}_{\mathcal{H}}=\pi_1^*(\mathcal{L}_{E})$ since $\pi_1$ is compatible with the pairing, therefore

\begin{equation*}
\begin{split}
\pi_{1!}(p_1^*(\P) \otimes  \mathcal{L}_{\mathcal{H}})&=\pi_{1!}p_1^*(\P) \otimes \mathcal{L}_E\\
&={p'_1}^*\pi_{1!}(\P) \otimes \mathcal{L}_E.
\end{split}
\end{equation*}

Thus $\pi_{1!}\Phi_{\mathcal{H}}(\P) = \Phi_{E}(\pi_{1!}(\P))$ as wanted. The Lemma is proved.

\vspace{.15in}

\noindent
\textit{End of proof of the Proposition}. We may now turn to the proof of the equivalence in the Proposition. Let $\P_{\Omega_1}$ be an object of $\mathcal{P}_{\d,\Omega_1}$. There exists a sequence $\i$ and a simple perverse sheaf $P_{\Omega_1}$ on $\mathcal{H}^{im}_{\i,\Omega_1}$ such that $\P=\pi_{1!}(P_{\Omega_1})$. Moreover, $P_{\Omega_1}$ is constant along the fibers of $e$, and is pulled back from $\mathcal{K}^{im}_{\i,\Omega_1}$, i.e.
$P_{\Omega_1}=e^*(P'_{\Omega_1})$ for some $P'_{\Omega_1}$. It follows that the restriction $j^*(\P_{\Omega_1})$ of $\P_{\Omega_1}$ to $E_{\d,\Omega_1}^{reg}$ is a simple perverse sheaf, and that moreover we have $j^*(\P_{\Omega_1})=\rho^*(P'_{\Omega_1})$. This may be summarized in the following diagram
$$\xymatrix{
P_{\Omega_1} \ar[r]^-{\pi_{1!}}& \mathbb{P}_{\Omega_1} \ar[d]^-{j^*}\\
P'_{\Omega_1} \ar[u]^-{e^*} \ar[r]_-{\rho^*} & j^*(\P_{\Omega_1})
}
$$
In particular, the restriction $j^*(\mathbb{P}_{\Omega_1})$ of $\P_{\Omega_1}$ completely determines the whole of $\P_{\Omega_1}$. The IC extension $j_{*!}(j^*(\P_{\Omega_1}))$ is a simple constituent of
$\P_{\Omega_1}$.

\vspace{.1in}

Now consider $\Phi_{E^{reg}}(j^*(\P_{\Omega_1}))$. This is a simple perverse sheaf on $E^{reg}_{\d,\Omega_2}$, and its IC extension belongs to $\overline{\mQ}_{\d,\Omega_2}$. By the discussion in the previous paragraph (applied to $\Omega_2$ instead of $\Omega_1$), it appears as the restriction $j^*(\mathbb{P}_{\Omega_2})$ of some (unique) object $\P_{\Omega_2}$ in $\mQ_{\d,\Omega_2}$. We will prove that $\Phi_E(\mathbb{P}_{\Omega_1}) \simeq \P_{\Omega_2}$. Indeed, there exists simple perverse sheaves $P_{\Omega_2}$ an $P'_{\Omega_2}$ on $\mathcal{H}^{im}_{\i,\Omega_2}$ and $\mathcal{K}^{im}_{\i,\Omega_2}$ satisfying $e^*(P'_{\Omega_2})=P_{\Omega_2}$ and $\rho^*(P'_{\Omega_2})=j^*(\P_{\Omega_2})$. Since $\rho^*: D^b(\mathcal{K}^{im}_{\i,\Omega_i}) \to D^b(E_{\d,\Omega_i}^{reg})$ is a fully faithful embedding for $i=1, 2$, we deduce from the Lemma that $\Phi_{\mathcal{K}}(P'_{\Omega_1})\simeq P'_{\Omega_2}$. Using the Lemma again, it follows that 
$$\Phi_{\mathcal{H}}(P_{\Omega_1})=\Phi_{\mathcal{H}}(e^*P'_{\Omega_1})=e^* \Phi_{\mathcal{K}}(P'_{\Omega_1})\simeq e^* P'_{\Omega_2}=P_{\Omega_2}.$$
But then 
$$\Phi_{E}(\mathbb{P}_{\Omega_1})=\Phi_E(\pi_{1!}P_{\Omega_1})=\pi_{1!} \Phi_{\mathcal{H}}(P_{\Omega_1}) \simeq \pi_{1!}P_{\Omega_2}=\mathbb{P}_{\Omega_2}$$
as wanted.

\vspace{.1in}

We have shown that $\Phi_E(\P)$ belongs to $\mQ_{\d,\Omega_2}$ for any object $\P$ of $\mQ_{\d,\Omega_1}$. The reverse inclusion is proved in a similar fashion, using the Fourier transform from $\Omega_2$ to $\Omega_1$. The Proposition is proved.

\vspace{.2in}

\paragraph{\textbf{3. Miscellaneous.}} i) The paternity of canonical basis theory is wrongly attributed in the introduction to our paper. We refer to \cite{EK} (introduction and end of section~2.2.) as well as to
\cite{Nag} for the correct history.

ii) In the proof of Proposition~5.1.  we introduce at some point a  
large commutative diagram containing four squares (after the definition of condition~$(\star)$ ). In the  paragraph which follows the diagram we state ``one deduces that the two squares in the above diagram are cartesian''. Instead, one should read : ``one deduces that the two leftmost squares in the above diagram are cartesian''.

\vspace{.2in}

\paragraph{\textbf{4. Relation to other works.}} In \cite{LusTight}, Lusztig gave a construction of an algebra $U$ equipped with a canonical basis $B$, attached to any quiver with loops. Our construction is closely related to his. Namely, he defines a category $\mathcal{Q}$ of \textit{simple} perverse sheaves on varieties of representations of $Q$, in a way similar to \cite{Lusbook} (as we did in \cite{KS}). He then defines $U$ as the (graded) Grothendieck group of $\mathcal{Q}$, and $B$ as the basis corresponding to the simple perverse sheaves. He doesn't impose the nilpotency condition along the loops as we do, but the two constructions are nicely related by a Fourier-Deligne transform along these loops. By definition, the algebra $U$ contains some natural elements $E_i^{(n)}$ for all vertices $i$ and positive integers $n$ : when $i$ is real then $E_i^{(n)}=E_i^{(1)}/[n]!$, but when $i$ is imaginary $E_i^{(n)}$ can not be expressed in terms of $E_i^{(1)}$. In particular, it is not known in general whether or not these elements generate $U$.

The algebra $\mathbf{U}$ which we construct is naturally a subalgebra of $U$ which is, by Theorem~4.1., generated by the elements $E_i^{(1)}$ for $i \in I$. Moreover, when $Q$ has no imaginary vertex with a single loop attached to it (such as the monster quiver for instance) then our canonical basis $\mathbf{B}$ consists of simple perverse sheaves and is thus a subset of $B$. In particular, we have $\mathbf{B}=B \cap \mathbf{U}$.

\vspace{.1in}

In the recent preprint \cite{LL}, Y. Li and Z. Lin further studied the canonical bases of GKM algebras defined in \cite{KS}. They showed that these bases are in a certain sense independent of the number of loops at imaginary vertices (provided there is at least one). In particular, the several possible choices for a quiver corresponding to a fixed GKM algebra all give rise to the same canonical bases. Some of the arguments used in the proof of the Proposition in the present paper also appear (at least implicitly) in \cite{LL}.

\vspace{.1in}

\centerline{\textbf{Acknowledgements.}}

\vspace{.1in}

The authors thank George Lusztig for his numerous remarks and comments. We are also grateful to
Michel van den Bergh for pointing out the reference \cite{LusTight} to us.

\vspace{.2in}

\small{

}

\vspace{.2in}

\end{document}